\documentclass[preprints,article,accept,moreauthors,pdftex]{my-mdpi} 

\firstpage{1} 
\pubvolume{5}
\issuenum{4}
\articlenumber{219}
\doinum{10.3390/fractalfract5040219}
\pubyear{2021}
\copyrightyear{2021}
\externaleditor{Academic Editor: Paul Eloe} 
\datereceived{31 October 2021} 
\daterevised{10 November 2021}
\dateaccepted{12 November 2021} 
\datepublished{14 November 2021} 
\hreflink{https://doi.org/10.3390/\linebreak fractalfract5040219} 

\pdfoutput=1

\usepackage{url}

\Title{Numerical Solution of Variable-Order Fractional Differential Equations Using Bernoulli Polynomials}

\TitleCitation{Numerical Solution of Variable-Order Fractional Differential Equations Using Bernoulli Polynomials}

\Author{Somayeh Nemati $^{1,\dagger}$\orcidA{}, 
Pedro M. Lima $^{2,\dagger}$\orcidB{} and 
Delfim F. M. Torres $^{3,}$*$^{,\dagger}$\orcidC{}}

\AuthorNames{Somayeh Nemati, Pedro M. Lima and Delfim F. M. Torres}

\AuthorCitation{Nemati, S.; Lima, P.M.; Torres, D.F.M.}

\address{$^{1}$ \quad Department of Applied Mathematics, Faculty of Mathematical Sciences, 
University of Mazandaran, P.O. Box 47416-95447 Babolsar, Iran; 
s.nemati@umz.ac.ir\\
$^{2}$ \quad Centro de Matem\'{a}tica Computacional e Estoc\'astica, 
Instituto Superior T\'{e}cnico, Universidade de Lisboa, 
Av.~Rovisco Pais, 1049-001 Lisboa, Portugal; 
pedro.t.lima@ist.utl.pt\\
$^{3}$ \quad Center for Research and Development in Mathematics and Applications (CIDMA), 
Department of Mathematics, University of Aveiro, 3810-193 Aveiro, Portugal}

\corres{Correspondence: delfim@ua.pt}

\firstnote{These authors contributed equally to this work.}

\abstract{We introduce a new numerical method, based on Bernoulli polynomials, 
for solving multiterm variable-order fractional differential equations. 
The variable-order fractional derivative was considered in the Caputo sense, 
while the Riemann--Liouville integral operator was used to give approximations 
for the unknown function and its variable-order derivatives. An operational matrix 
of variable-order fractional integration was introduced for the Bernoulli functions.
By assuming that the solution of the problem is sufficiently smooth, 
we approximated a given order of its derivative using Bernoulli polynomials. 
Then, we used the introduced operational matrix to find some approximations 
for the unknown function and its derivatives. Using these approximations 
and some collocation points, the problem was reduced to the solution of a system 
of nonlinear algebraic equations. An error estimate is given for the approximate 
solution obtained by the proposed method. Finally, 
five illustrative examples were considered to demonstrate 
the applicability and high accuracy of the proposed technique,
comparing our results with the ones obtained 
by existing methods in the literature and
making clear the novelty of the work. The numerical results
showed that the new method is efficient, giving high-accuracy 
approximate solutions even with a small number of basis functions 
and when the solution to the problem is not infinitely differentiable,
providing better results and a smaller number of basis 
functions when compared to state-of-the-art methods.}

\keyword{fractional differential equations;
numerical methods;
variable-order fractional calculus;
operational matrix of variable-order fractional integration;
Bernoulli polynomials}

\MSC{34A08; 65L60}

\begin{document}


\section{Introduction}
\label{sec:1}

In the last few decades, fractional calculus has attracted the attention 
of many scientists in different fields such as mathematics, physics, 
chemistry, and engineering. Due to the fact that fractional operators 
consider the evolution of the system, by taking the global correlation, 
and not only local characteristics, some physical phenomena are better 
described by fractional derivatives~\cite{Almeida}. There are many definitions 
of fractional differentiation and integration in the literature 
(for details, see, e.g.,~\cite{Samko,Podlubny,Kilbas,Mainardi}). 
However, the most commonly used are the definitions of the Caputo derivative 
and Riemann--Liouville fractional integral operators. For 
new fractional derivatives with nonlocal and nonsingular
kernels, with applications in rheological models, 
we refer the reader to~\cite{rev01:01}.

A recent generalization of the theory of fractional calculus is to allow 
the fractional order of the derivatives to be dependent on time, i.e., 
to be nonconstant or of variable order. In~\cite{Samko1}, the authors 
investigated operators when the order of the fractional derivative is 
variable with time~\cite{Odzijewicz,Chen}. 

The nonlocal properties of systems are more visible 
with variable-order fractional calculus, 
and many real-world phenomena in physics, mechanics, control, 
and signal processing have been described by this approach 
\cite{Coimbra,Odzijewicz1,Ostalczyk,Rapaic}. In particular,
several applications of variable-order fractional calculus 
are found in engineering mechanics; see~\cite{rev03:06}
for an application of variable-order fractional operators 
to model the microscopic structure of a material,~\cite{rev03:07} 
for an application of the Riesz--Caputo
fractional derivative of space-dependent order in
continuum elasticity,~\cite{rev03:08,rev03:10} for the nonlinear 
viscoelastic behavior of fractional systems with variable 
time-dependent fractional order, and~\cite{rev03:09}
for the use of variable-order fractional calculus
in the static response of nonlocal beams having 
either a porous or a functionally graded core.

Obtaining analytic solutions for such fractional differential equations (FDEs) 
is, however, very difficult. Therefore, in most cases, the exact solution is not known, 
and one needs to seek a numerical approximation. Therefore, many 
researchers have introduced and developed numerical methods in order 
to obtain approximated solutions for this class of equations. For example,
in~\cite{El-Sayed,Wang}, Legendre polynomials were used to construct 
a numerical solution for a class of multiterm variable-order FDEs,
while~\cite{Chen2} used Legendre wavelets and operational matrices;
References~\cite{Liu,Nagy} used Chebyshev polynomials, respectively of the second 
and fourth kinds; Reference~\cite{Chen1} used Bernstein polynomials.
The book~\cite{Tavares} showed the usefulness of numerical methods
for approximating variable-order fractional operators in the framework 
of the calculus of variations; the paper~\cite{Shen} adopted Coimbra's 
variable-order time-fractional operator, discussing the stability, convergence 
and solvability of a numerical scheme based on Fourier analysis.
Reference~\cite{rev2:01} implemented a numerical method for solving 
a circulant Halvorsen system described by Caputo fractional variable-order derivatives.
For multivariable-order differential equations with nonlocal and 
nonsingular kernels, we refer to~\cite{rev01:03}, where a collocation 
method was developed based on Chebyshev polynomials of the fifth kind. 

Here, we considered the following general form 
of a multiterm variable-order FDE:
\begin{equation}
\label{1.1}
{_0^CD}_t^{\alpha(t)}y(t)=F\left(t,y(t),{_0^CD_t^{\alpha_1(t)}}y(t),
{_0^CD_t^{\alpha_2(t)}} y(t),\ldots,{_0^CD_t^{\alpha_k(t)}}y(t)\right),
\quad 0<t\leq 1,
\end{equation}
with initial conditions
\begin{equation}
\label{1.2}
y^{(i)}(0)=y_0^i, \quad i=0,1,\ldots,n-1,
\end{equation}
where $n$ is the smallest positive integer number such that
for all $t\in[0,1]$, one has $0<\alpha(t)\leq n$, 
$0<\alpha_1(t)<\alpha_2(t)<\ldots <\alpha_k(t)<\alpha(t)$, 
and ${_0^CD_t^{\alpha(t)}}$, ${_0^CD_t^{\alpha_1(t)}}$,\ldots, ${_0^CD_t^{\alpha_k(t)}}$ 
are the (left) fractional derivatives of variable-order defined in the Caputo sense. 
Problems of the form \eqref{1.1} and \eqref{1.2} have a practical impact. 
Indeed, specific applications are found in noise reduction and signal processing 
\cite{rev03:01,rev03:02}, the processing of geographical data~\cite{rev03:03}, 
and signature verification~\cite{rev03:04}.

In recent years, Bernoulli polynomials have been shown to be a powerful mathematical 
tool in dealing with various problems of a dynamical nature, e.g.,
for solving numerically high-order Fredholm integrodifferential equations
\cite{Bhrawy}, pantograph equations~\cite{Tohidi},
partial differential equations~\cite{Toutounian},
linear Volterra and nonlinear Volterra--Fredholm--Hammerstein integral equations 
\cite{Bazm}, as well as optimal control problems~\cite{Keshavarz}. Here,
we employ a spectral method based on Bernoulli polynomials 
in order to obtain numerical solutions to the problem 
\eqref{1.1} and \eqref{1.2}. Our method consists of reducing the problem to 
a system of nonlinear algebraic equations. To do this, we introduced an accurate 
operational matrix of variable-order fractional integration for 
the Bernoulli polynomials' basis vector. To the best of our knowledge, 
this is the first time in the literature that such a method 
for solving a general class of multiterm variable-order FDEs based on the 
Riemann--Liouville fractional integral of the basis vector has been introduced. 

The rest of this paper is organized as follows. In Section~\ref{sec:2}, 
some preliminaries of variable-order fractional calculus are reviewed 
and some properties of the Bernoulli polynomials are recalled. 
Section~\ref{sec:3} is devoted to introducing the operational matrix 
of variable-order fractional integration for Bernoulli polynomials. 
In Section~\ref{sec:4}, we present a new numerical method for 
solving the problem \eqref{1.1} and \eqref{1.2} by using the operational matrix 
technique and collocation points. Section~\ref{sec:5} is concerned with 
presenting an error estimate for the numerical solution obtained by this 
new scheme. In Section~\ref{sec:6}, several multiterm variable-order 
FDEs are considered and solved, using the introduced method. Finally,
concluding remarks are given in Section~\ref{sec:7}, where some possible 
future directions of research are also pointed out. 
 

\section{Preliminaries}
\label{sec:2}

In this section, a brief review of the necessary definitions 
and properties of the variable-order fractional calculus is presented. 
Furthermore, Bernoulli polynomials, and some of their properties, are recalled.


\subsection{Some Preliminaries of Variable-Order Fractional Calculus}

We followed the notations of~\cite{Almeida}.

\begin{Definition}[See, e.g.,~\cite{Almeida}]
The left Riemann--Liouville fractional integral of order $\alpha(t)$ is defined~by: 
\begin{equation*}
_0I_t^{\alpha(t)}y(t)=\frac{1}{\Gamma{(\alpha(t))}}
\int_0^t(t-s)^{\alpha(t)-1}y(s)ds,
\quad t>0,
\end{equation*}
where $\Gamma(\cdot)$ is the Euler gamma function.
\end{Definition}

\begin{Lemma}[See Chapter~1, Lemma~8, of~\cite{Almeida}]
\label{lemma1}
Let $y$ be the power function $y(t)=t^\nu$. Then, for $\nu>-1$, we have:
\begin{equation*}
_0I_t^{\alpha(t)}y(t)
=\frac{\Gamma(\nu+1)}{\Gamma{(\nu+1+\alpha(t))}}
t^{\nu+\alpha(t)},\quad \alpha(t)\geq 0.
\end{equation*}
\end{Lemma}

\begin{Definition}[See, e.g.,~\cite{Almeida}]
The left Caputo fractional derivative of order $\alpha(t)$ is defined by: 
\begin{equation*}
\begin{split}
&{_0^CD}_t^{\alpha(t)}y(t)=\frac{1}{\Gamma(n-\alpha(t))}
\int_0^t(t-s)^{n-\alpha(t)-1}y^{(n)}(s)ds,\quad n-1<\alpha(t)<n,\\
&{_0^CD}_t^{\alpha(t)}y(t)=y^{(n)}(t),\quad \alpha(t)=n.
\end{split}
\end{equation*}
\end{Definition}

For $0\leq\alpha(t)\leq n$, $n\in\mathbb{N}$, and $\gamma>0$, 
two useful properties of the Caputo derivative 
and Riemann--Liouville integral are:
\begin{equation}
\label{2.1}
_0I_t^{\gamma}({_0^CD}_t^{\gamma}y(t))
=y(t)-\sum_{i=0}^{\lceil\gamma\rceil-1}y^{(i)}(0)\frac{t^i}{i!},
\quad t>0,
\end{equation}
\begin{equation}
\label{2.2}
_0I_t^{n-\alpha(t)}(y^{(n)}(t))
={_0^CD}_t^{\alpha(t)}y(t)
-\sum_{i=\lceil\alpha(t)\rceil}^{n-1}y^{(i)}(0)
\frac{t^{i-\alpha(t)}}{\Gamma(i+1-\alpha(t))},
\quad t>0,
\end{equation}
where $\lceil \cdot \rceil$ is the ceiling function.


\subsection{Bernoulli Polynomials}

Bernoulli polynomials build a family of independent polynomials 
that form a complete basis for the space $L^2[0,1]$, which is 
the space of all square integrable functions on the interval $[0,1]$. 
The Bernoulli polynomial of degree $m$, $\beta_m(t)$, 
is defined as follows~\cite{Costabile}:
\begin{equation}
\label{2.3}
\beta_m(t)=\sum_{i=0}^m \binom{m}{i}b_{m-i}t^i,
\end{equation}
where $b_k$, $k=0,1,\ldots,m$, are the Bernoulli numbers 
that appear in the series expansion of trigonometric functions~\cite{Arfken} 
and can be defined by the following identity: 
\vspace{6pt}
\begin{equation*}
\frac{t}{e^t-1}=\sum_{i=0}^\infty b_i \frac{t^i}{i!}.
\end{equation*}

The first four Bernoulli polynomials are:
\begin{equation*}
\begin{split}
&\beta_0(t)=1,\\
&\beta_1(t)=t-\frac{1}{2},\\
&\beta_2(t)=t^2-t+\frac{1}{6},\\
&\beta_3(t)=t^3-\frac{3}{2}t^2+\frac{1}{2}t.
\end{split}
\end{equation*}

The Bernoulli polynomials satisfy the following property~\cite{Arfken}:
\begin{equation}
\label{2.4}
\int_0^1\beta_i(t)\beta_j(t)dt
=(-1)^{i-1}\frac{i!j!}{(i+j)!}b_{i+j},\quad i,j\geq 1.
\end{equation}

Any arbitrary function $y\in L^2[0,1]$ can be approximated 
using the Bernoulli polynomials as
\begin{equation}
\label{2.5}
y(t)\simeq \sum_{m=0}^M a_m\beta_m(t)=A^TB(t),
\end{equation}
where
\begin{equation}
\label{2.6}
B(t)=[\beta_0(t),\beta_1(t),\ldots,\beta_M(t)]^T
\end{equation}
and
\begin{equation*}
A=[a_0,a_1,\ldots,a_M]^T.
\end{equation*}

The coefficient vector $A$ in \eqref{2.5} is calculated 
by the following formula (see~\cite{Keshavarz}):
\begin{equation}
\label{2.7}
A=D^{-1}\langle y(t),B(t)\rangle,
\end{equation}
where $\langle \cdot,\cdot \rangle$ is the inner product, 
defined for two arbitrary functions $f,g\in L^2 [0,1]$ as
\begin{equation*}
\langle f(t),g(t) \rangle=\int_0^1f(t)g(t)dt
\end{equation*}
and $D=\langle B(t),B(t)\rangle$, 
which is calculated using \eqref{2.4}.


\section{Operational Matrix of Variable-Order Fractional Integration}
\label{sec:3}

The aim of this section is to introduce an accurate operational matrix 
of variable-order fractional integration for Bernoulli functions. 
To do this, we display the Bernoulli basis vector $B(t)$, 
given by \eqref{2.6}, in terms of the Taylor basis functions, as follows:
\begin{equation}
\label{3.1}
B(t)=Q\mathbb{T}(t),
\end{equation}
where $\mathbb{T}$ is the Taylor basis vector
\begin{equation*}
\mathbb{T}(t)=[1,t,t^2,\ldots,t^M]^T,
\end{equation*}
and $Q$ is the change-of-basis matrix, 
which is obtained using \eqref{2.3} as
\vspace{6pt}
\begin{equation*}
Q=\left[\begin{array}{ccccccc}
1& 0 & 0 & 0 & 0 &\ldots & 0\\
-\frac{1}{2} & 1 & 0 & 0 & 0&\ldots & 0\\
\frac{1}{6} & -1 & 1 & 0 & 0&\ldots & 0\\
0 &\frac{1}{2}& -\frac{3}{2} & 1 &0 &\ldots & 0\\
\vdots & \vdots& \vdots & \vdots &\vdots& & \vdots\\
b_M & \binom{M}{1}b_{M-1} & \binom{M}{2}b_{M-2} 
& \binom{M}{3}b_{M-3} & \binom{M}{4}b_{M-4}&\ldots & 1\\
\end{array}
\right].
\end{equation*}

It is obvious that the matrix $Q$ is a nonsingular matrix. 
Therefore, we can write:
\begin{equation}
\label{3.2}
\mathbb{T}(t)=Q^{-1}B(t).
\end{equation}

Taking \eqref{3.1} into account and by applying the left 
Riemann--Liouville fractional integral operator of order 
$\alpha(t)$ to the vector $B(t)$, we obtain
\begin{equation}
\label{3.3}
_0I_t^{\alpha(t)} B(t)={_0I_t^{\alpha(t)}}(Q\mathbb{T}(t))
=Q ({_0I_t^{\alpha(t)}}\mathbb{T}(t))=QS_t^{\alpha(t)}\mathbb{T}(t),
\end{equation}
where $S_t^{\alpha(t)}$ is a diagonal matrix, 
which is obtained using Lemma~\ref{lemma1} as follows:
\begin{equation*}
S_t^{\alpha(t)}=\left[\begin{array}{cccccc}
\frac{1}{\Gamma(1+\alpha(t))}t^{\alpha(t)} & 0 & 0 &0&\cdots &0\\
0 & \frac{1}{\Gamma(2+\alpha(t))}t^{\alpha(t)} & 0 &0&\cdots &0\\
0 & 0 & \frac{2}{\Gamma(3+\alpha(t))}t^{\alpha(t)} &0&\cdots &0\\
\vdots & \vdots & \vdots & \vdots& & \vdots\\
0 & 0 & 0 & 0 & \cdots & \frac{\Gamma(M+1)}{\Gamma{(M+1+\alpha(t))}}t^{\alpha(t)}\\
\end{array}
\right].
\end{equation*}

Finally, by using \eqref{3.2} into \eqref{3.3}, we have: 
\begin{equation}
\label{3.4}
_0I_t^{\alpha(t)} B(t)=QS_t^{\alpha(t)}Q^{-1}B(t)
=P_t^{\alpha(t)}B(t),
\end{equation}
where $P_t^{\alpha(t)}=QS_t^{\alpha(t)}Q^{-1}$ 
is a matrix of dimension $(M+1)\times (M+1)$, which we call
the operational matrix of variable-order fractional integration 
of order $\alpha(t)$ of Bernoulli functions. Since $Q$ and $Q^{-1}$ 
are lower triangular matrices and $S_t^{\alpha(t)}$ is a diagonal matrix, 
it is obvious that $P_t^{\alpha(t)}$ is also a lower triangular matrix. 
For example, with $M=2$, one has:
\begin{equation*}
P_t^{\alpha(t)}=\left[
\begin{array}{ccc}
p_{1,1} & 0 & 0 \\
p_{2,1} & \frac{1}{\Gamma (\alpha (t)+2)}t^{\alpha (t)} & 0 \\
p_{3,1} & \left(\frac{2}{\Gamma (\alpha (t)+3)}-\frac{1}{\Gamma (\alpha (t)+2)} 
\right)t^{\alpha (t)} & \frac{2}{\Gamma (\alpha (t)+3)} t^{\alpha (t)} \\
\end{array}
\right],
\end{equation*}
where
\begin{gather*}
p_{1,1} = \frac{1}{\Gamma (\alpha (t)+1)}t^{\alpha (t)},\\
p_{2,1} = \left(\frac{1}{2 \Gamma (\alpha (t)+2)}
-\frac{1}{2 \Gamma (\alpha (t)+1)}\right) t^{\alpha (t)},\\
p_{3,1} = \left(\frac{1}{6 \Gamma (\alpha (t)+1)}
-\frac{1}{2 \Gamma (\alpha (t)+2)}+\frac{2 }{3 
\Gamma (\alpha (t)+3)}\right)t^{\alpha (t)}.
\end{gather*}


\section{Numerical Method}
\label{sec:4}

This section is devoted to presenting a new numerical method 
for solving the multiterm variable-order FDE \eqref{1.1} 
with initial conditions \eqref{1.2}. To this aim, we set:
\begin{equation*}
n :=\max_{0< t\leq 1}\{\lceil \alpha(t)\rceil\}.
\end{equation*}

Then, by assuming $y\in C^n[0,1]$, we considered an approximation 
of the $n$-th order derivative of the unknown function $y$ 
using the Bernoulli functions as follows:
\begin{equation}
\label{4.1}
y^{(n)}(t)=A^TB(t),
\end{equation}
where $A$ is an $ (M+1)\times 1$ vector with unknown elements and $B(t)$ 
is the Bernoulli basis vector given by \eqref{2.6}. Taking into account 
the initial conditions \eqref{1.2} and using \eqref{2.1}, \eqref{3.4} 
and \eqref{4.1}, we obtain:
\begin{equation}
\label{4.2}
\begin{split}
y(t)&={_0I_t^{n}}(y^{(n)}(t))+\sum_{i=0}^{n-1}y^{(i)}(0)\frac{t^i}{i!}\\
&=A^T({_0I_t^{n}}B(t))+\sum_{i=0}^{n-1}y_0^{i}\frac{t^i}{i!}\\
&=A^TP_t^{n}B(t)+\sum_{i=0}^{n-1}y_0^{i}\frac{t^i}{i!}.
\end{split}
\end{equation}

In a similar way, using \eqref{2.2}, \eqref{3.4}, and \eqref{4.1}, we obtain
\begin{equation}
\label{4.3}
{_0^CD}_t^{\alpha(t)}y(t)=A^TP_t^{n-\alpha(t)}B(t)
+\sum_{i=\lceil\alpha(t)\rceil}^{n-1}y_0^{i}
\frac{t^{i-\alpha(t)}}{\Gamma(i+1-\alpha(t))}
\end{equation}
and
\begin{equation}
\label{4.4}
{_0^CD}_t^{\alpha_j(t)}y(t)=A^TP_t^{n-\alpha_j(t)}B(t)
+\sum_{i=\lceil\alpha_j(t)\rceil}^{n-1}y_0^{i}
\frac{t^{i-\alpha_j(t)}}{\Gamma(i+1-\alpha_j(t))},
\quad j=1,\ldots,k.
\end{equation}

By substituting the approximations given in \eqref{4.2}--\eqref{4.4} 
into Equation \eqref{1.1}, we have:
\begin{multline}
\label{4.5}
F\left(t,A^TP_t^{n}B(t)+\sum_{i=0}^{n-1}y_0^{i}
\frac{t^i}{i!},A^TP_t^{n-\alpha_1(t)}B(t)
+\sum_{i=\lceil\alpha_1(t)\rceil}^{n-1}y_0^{i}
\frac{t^{i-\alpha_1(t)}}{\Gamma(i+1-\alpha_1(t))},\right.\\
\left. \ldots, A^TP_t^{n-\alpha_k(t)}B(t)
+\sum_{i=\lceil\alpha_k(t)\rceil}^{n-1}y_0^{i}
\frac{t^{i-\alpha_k(t)}}{\Gamma(i+1-\alpha_k(t))}\right)\\
= A^TP_t^{n-\alpha(t)}B(t)+\sum_{i=\lceil\alpha(t)\rceil}^{n-1}
y_0^{i}\frac{t^{i-\alpha(t)}}{\Gamma(i+1-\alpha(t))}.
\end{multline}

By using \eqref{4.5} at $M+1$ collocation points, 
which are chosen as $t_j=\frac{j+1}{M+2}$, with \linebreak $j=0,1,\ldots,M$, 
we obtain the following system of nonlinear algebraic equations:
\begin{equation}
\label{4.6}
\begin{split}
A^TP_{t_j}^{n-\alpha({t_j})}B({t_j})
&+\sum_{i=\lceil\alpha({t_j})\rceil}^{n-1}
y_0^{i}\frac{{t_j}^{i-\alpha({t_j})}}{\Gamma(i+1-\alpha({t_j}))}\\
&=F\left({t_j},A^TP_{t_j}^{n}B({t_j})+\sum_{i=0}^{n-1}y_0^{i}
\frac{{t_j}^i}{i!},A^TP_{t_j}^{n-\alpha_1({t_j})}B({t_j})\right.\\
&\quad +\sum_{i=\lceil\alpha_1({t_j})\rceil}^{n-1}y_0^{i}
\frac{{t_j}^{i-\alpha_1({t_j})}}{\Gamma(i+1-\alpha_1({t_j}))},
\ldots,A^TP_{t_j}^{n-\alpha_k({t_j})}B({t_j})\\
&\quad \left.+\sum_{i=\lceil\alpha_k({t_j})\rceil}^{n-1}y_0^{i}
\frac{{t_j}^{i-\alpha_k({t_j})}}{\Gamma(i+1-\alpha_k({t_j}))}\right).
\end{split}
\end{equation}

System \eqref{4.6} includes $M+1$ nonlinear algebraic equations 
in terms of the unknown parameters of vector $A$. After solving 
this system, an approximation of the solution of the problem 
\eqref{1.1} and \eqref{1.2} is given by \eqref{4.2}.


\section{Error Estimate}
\label{sec:5}

The purpose of this section is to obtain an estimate of the error norm 
for the approximate solution obtained by the proposed method in 
Section~\ref{sec:4}. We assume that $f(t)=y^{(n)}(t)$ is a sufficiently 
smooth function on $[0,1]$ and $q_M(t)$ is the interpolating polynomial 
to $f$ at points $t_s$, where $t_s$, $s=0,1,\ldots,M$, are the roots 
of the $(M + 1)$-degree shifted Chebyshev polynomial in $[0,1]$. 
Then, according to the interpolation error, we have: 
\begin{equation*}
f(t)-q_M(t)=\frac{f^{(M+1)}(\xi)}{(M+1)!}\prod_{s=0}^M(t-t_s),\quad \xi\in(0,1).
\end{equation*}

Therefore,
\begin{equation}
\label{5.2}
|f(t)-q_M(t)|\leq\max_{t\in (0,1)}
\frac{\left|f^{(M+1)}(t)\right|}{(M+1)!}\prod_{s=0}^M|t-t_s|.
\end{equation}

We assume that there is a real number $\kappa$ such that:
\begin{equation}
\label{5.3}
\max_{t\in (0,1)}\left|f^{(M+1)}(t)\right|\leq \kappa.
\end{equation}

By using \eqref{5.3} in \eqref{5.2} and taking into consideration 
the estimates for Chebyshev interpolation nodes~\cite{Mason}, we obtain: 
\begin{equation}
\label{5.4}
|f(t)-q_M(t)|\leq\frac{\kappa}{2^{2M+1}(M+1)!}.
\end{equation}

As a consequence of \eqref{5.4}, we obtain the following result.

\begin{Theorem}
\label{th1}
Let $f_M(t)=A^TB(t)$ be the Bernoulli functions' expansion of 
a sufficiently smooth function $f$ defined on $[0,1]$, 
where $A$ and $B(t)$ are given, respectively, by \eqref{2.6} 
and \eqref{2.7}. Then, there exists a real number $\kappa$ such that:
\begin{equation}
\label{5.5}
\|f(t)-f_M(t)\|\leq \frac{\kappa}{2^{2M+1}(M+1)!}.
\end{equation}
\end{Theorem}

\begin{proof}
Let $\Pi_M$ be the space of all polynomials of degree $\leq M$ on $t\in[0,1]$. 
By definition, $f_M$ is the best approximation of $f$ in $\Pi_M$. Therefore, we have:
\begin{equation*}
\|f(t)-f_M(t)\|_2\leq \|f(t)-g(t)\|_2,
\end{equation*}
where $g$ is any arbitrary polynomial in $\Pi_M$. Therefore, we can write:
\begin{equation}
\label{5.6}
\|f(t)-f_M(t)\|_2^2=\int_0^1|f(t)-f_M(t)|^2dt\leq\int_0^1|f(t)-q_M(t)|^2dt,
\end{equation}
where $q_M$ is the interpolating polynomial of $f$, as discussed before. 
As a result, taking into consideration \eqref{5.4} in \eqref{5.6}, 
we easily obtain \eqref{5.5}.
\end{proof}

With the help of Theorem~\ref{th1} and by assuming $f(t)=y^{(n)}(t)$, 
we obtain the following result.

\begin{Theorem}
\label{th2}
Let the exact solution $y$ of the problem \eqref{1.1} and \eqref{1.2} 
be a real-valued sufficiently smooth function and $y_M$ be the 
approximate solution of this problem obtained by the method 
proposed in Section~\ref{sec:4}. Then, we have:
\begin{equation*}
\|y(t)-y_M(t)\|_2\leq \frac{\kappa}{2^{2M+1}(M+1)!(n-1)!\sqrt{2n(2n-1)}},
\end{equation*}
where $y$ is the exact solution and 
$\kappa=\max_{t\in (0,1)}\left|y^{(nM+n)}(t)\right|$, 
$n=\max_{0< t\leq 1}\{\lceil \alpha(t)\rceil\}$.
\end{Theorem}

\begin{proof}
Let $X$ be the space of all real-valued functions defined on $[0,1]$ 
and $_0I_t^{n}:X\rightarrow X$ be the Riemann--Liouville integral operator. 
We use the following definition for the norm of the operator $_0I_t^{n}$:
\begin{equation*}
\|_0I_t^{n}\|_2=\sup_{\|g\|_2=1}\|_0I_t^{n}g\|_2.
\end{equation*}

To continue the proof, first, we introduce an upper bound for 
$\|_0I_t^{n}\|_2$. To this end, using the definition of the left 
Riemann--Liouville integral operator and Schwarz's inequality, we~obtain:
\begin{equation*}
\begin{split}
\left\|_0I_t^{n}g\right\|_2^2
&=\left\|\frac{1}{(n-1)!}\int_0^t(t-s)^{n-1}g(s)ds\right\|_2^2\\
&=\frac{1}{\left[(n-1)!\right]^2}\left\|\int_0^t(t-s)^{n-1}g(s)ds\right\|_2^2
=\frac{1}{\left[(n-1)!\right]^2}
\int_0^1\left|\int_0^t(t-s)^{n-1}g(s)ds\right|^2dt\\
&\leq \frac{1}{\left[(n-1)!\right]^2}\int_0^1\left(
\int_0^t(t-s)^{2n-2}ds\right)\left(\int_0^1|g(s)|^2ds\right)dt\\
&= \frac{1}{\left[(n-1)!\right]^2 2n(2n-1)}.
\end{split}
\end{equation*}

Therefore, we have: 
\begin{equation}
\label{5.7}
\left\|_0I_t^{n}\right\|_2\leq \frac{1}{(n-1)!\sqrt{2n(2n-1)}}.
\end{equation}

On the other hand, from Theorem~\ref{th1}, we have the following error bound:
\begin{equation}
\label{5.8}
\left\|y^{(n)}(t)-A^TB(t)\right\|_2\leq \frac{\kappa}{2^{2M+1}(M+1)!}.
\end{equation}

Finally, using \eqref{2.1}, \eqref{4.2}, \eqref{5.7} and \eqref{5.8}, we obtain that:
\begin{equation*}
\begin{split}
\left\|y(t)-y_M(t)\right\|_2
&=\left\|{_0I_t^{n}}(y^{(n)}(t))+\sum_{i=0}^{n-1}y_0^{i}
\frac{t^i}{i!}-\left( {_0I_t^{n}}(A^TB(t))
+\sum_{i=0}^{n-1}y_0^{i}\frac{t^i}{i!}\right)\right\|_2\\
&=\left\|{_0I_t^{n}}(y^{(n)}(t)-A^TB(t))\right\|_2\\
&\leq \left\|{_0I_t^{n}} \right\|_2\left\|y^{(n)}(t)-A^TB(t)\right\|_2\\
&\leq \frac{\kappa}{2^{2M+1}(M+1)!(n-1)!\sqrt{2n(2n-1)}},
\end{split}
\end{equation*}
which completes the proof.
\end{proof}


\section{Illustrative Examples}
\label{sec:6}

In this section, we apply our method to some variable-order FDEs and, moreover, 
to one variable-order fractional pantograph differential equation 
(Example~\ref{ex4}), comparing the results with the ones obtained 
by existing methods in the literature. We implemented our method 
and performed our numerical simulations with \textsf{Mathematica 12}. 


\begin{Example}
\label{ex1}
\textls[-15]{In our first example, we considered the following 
multiterm variable-order FDE~\cite{El-Sayed,Nagy,Liu}:}
\begin{equation}
\label{6.1}
{_0^CD}_t^{2t}y(t)+t^{\frac{1}{2}}{_0^CD}_t^{\frac{t}{3}}y(t)
+t^{\frac{1}{3}}{_0^CD}_t^{\frac{t}{4}}y(t)
+t^{\frac{1}{4}}{_0^CD}_t^{\frac{t}{5}}y(t)
+t^{\frac{1}{5}}y(t)=g(t),
\quad 0<t\leq 1,
\end{equation}
where
\[
g(t)=-\frac{t^{2-2t}}{\Gamma(3-2t)}-t^{\frac{1}{2}}
\frac{ t^{2-\frac{t}{3}}}{\Gamma (3-\frac{t}{3})}
-t^{\frac{1}{3}}\frac{t^{2-\frac{t}{4}}}{\Gamma (3-\frac{t}{4})}
-t^{\frac{1}{4}}\frac{t^{2-\frac{t}{5}}}{\Gamma (3-\frac{t}{5})}
+t^{\frac{1}{5}}\left(2-\frac{t^2}{2}\right)
\]
with initial conditions $y(0)=2$ and $y'(0)=0$. 
The exact solution to this problem is $y(t)=2-\frac{t^2}{2}$. 
As can be seen, we have $\alpha(t)=2t$. Therefore, 
to implement the proposed method, we introduce:
\begin{equation*}
n=\max_{0< t\leq 1}\{\lceil 2t\rceil\}=2.
\end{equation*}

We set $M=1$ and suppose
\begin{equation*}
y{''}(t)=A^TB(t),
\end{equation*}
where
\begin{equation*}
A=\left[a_0,a_1\right]^T \text{ and }B(t)=\left[1,t-\frac{1}{2}\right]^T.
\end{equation*}

The operational matrices of variable-order fractional integration are given as follows:
\begin{equation*}
P_t^{2-2t}=\left[
\begin{array}{cc}
\frac{t^{2-2 t}}{\Gamma (3-2 t)} & 0 \\
\frac{t^{2-2 t}}{2 \Gamma (4-2 t)}
-\frac{t^{2-2 t}}{2 \Gamma (3-2 t)} & \frac{t^{2-2 t}}{\Gamma (4-2 t)}
\end{array}
\right],\quad P_t^{2-\frac{t}{3}}=\left[
\begin{array}{cc}
\frac{t^{2-\frac{t}{3}}}{\Gamma \left(3-\frac{t}{3}\right)} & 0 \\
\frac{t^{2-\frac{t}{3}}}{2 \Gamma \left(4-\frac{t}{3}\right)}
-\frac{t^{2-\frac{t}{3}}}{2 \Gamma \left(3-\frac{t}{3}\right)} 
& \frac{t^{2-\frac{t}{3}}}{\Gamma \left(4-\frac{t}{3}\right)}
\end{array}
\right],
\end{equation*}
\begin{equation*}
P_t^{2-\frac{t}{4}}=\left[
\begin{array}{cc}
\frac{t^{2-\frac{t}{4}}}{\Gamma \left(3-\frac{t}{4}\right)} & 0 \\
\frac{t^{2-\frac{t}{4}}}{2 \Gamma \left(4-\frac{t}{4}\right)}
-\frac{t^{2-\frac{t}{4}}}{2 \Gamma \left(3-\frac{t}{4}\right)} 
& \frac{t^{2-\frac{t}{4}}}{\Gamma \left(4-\frac{t}{4}\right)}
\end{array}
\right],\quad P_t^{2-\frac{t}{5}}=\left[
\begin{array}{cc}
\frac{t^{2-\frac{t}{5}}}{\Gamma \left(3-\frac{t}{5}\right)} & 0 \\
\frac{t^{2-\frac{t}{5}}}{2 \Gamma \left(4-\frac{t}{5}\right)}
-\frac{t^{2-\frac{t}{5}}}{2 \Gamma \left(3-\frac{t}{5}\right)} 
& \frac{t^{2-\frac{t}{5}}}{\Gamma \left(4-\frac{t}{5}\right)}
\end{array}
\right].
\end{equation*}

Furthermore, we have: 
\begin{equation*}
P_t^{2}=\left[
\begin{array}{cc}
\frac{t^2}{2} & 0 \\
-\frac{t^2}{6} & \frac{t^2}{6}
\end{array}
\right].
\end{equation*}

Now, using the initial conditions and the aforementioned operational matrices, 
we obtain the following approximations for $y(t)$ and its variable-order derivatives:
\begin{equation*}
y(t)=A^TP_t^{2}B(t)+2,\quad {_0^CD}_t^{2t}y(t)=A^TP_t^{2-2t}B(t)
+\sum _{i=\lceil 2t\rceil }^{1} y_0^i\frac{ t^{i-2t}}{\Gamma (i-2t+1)},
\end{equation*}
\vspace{-10pt}
\begin{gather*}
{_0^CD}_t^{\frac{t}{3}}y(t)=A^TP_t^{2-\frac{t}{3}}B(t)
+\sum _{i=\lceil \frac{t}{3}\rceil }^{1} y_0^i
\frac{ t^{i-\frac{t}{3}}}{\Gamma (i-\frac{t}{3}+1)},\\
{_0^CD}_t^{\frac{t}{4}}y(t)=A^TP_t^{2-\frac{t}{4}}B(t)
+\sum _{i=\lceil \frac{t}{4}\rceil }^{1} y_0^i
\frac{ t^{i-\frac{t}{4}}}{\Gamma (i-\frac{t}{4}+1)},\\
{_0^CD}_t^{\frac{t}{5}}y(t)=A^TP_t^{2-\frac{t}{5}}B(t)
+\sum _{i=\lceil \frac{t}{5}\rceil }^{1} y_0^i
\frac{ t^{i-\frac{t}{5}}}{\Gamma (i-\frac{t}{5}+1)}.
\end{gather*}

By substituting these approximations into \eqref{6.1}, collocating 
the resulting equation at $t_0=\frac{1}{3}$, $t_1=\frac{2}{3}$, 
and finally, solving the obtained system, we obtain:
\begin{equation*}
a_0=-1,\quad a_1=0.
\end{equation*}

Therefore, we have:
\begin{equation*}
y(t)=[-1 \quad 0]\left[
\begin{array}{cc}
\frac{t^2}{2} & 0 \\
-\frac{t^2}{6} & \frac{t^2}{6} 
\end{array}
\right]\left[
\begin{array}{c}
1 \\
t-\frac{1}{2}
\end{array}
\right]+2=2-\frac{t^2}{2},
\end{equation*}
which is the exact solution. In this case, since the exact solution 
is a second-order polynomial, we can obtain it by applying 
the numerical method with just two basis functions. 
\end{Example}


\begin{Example}
\label{ex2} 
In our second example, we considered the following nonlinear 
variable-order FDE borrowed from~\cite{Hassan}:
\begin{equation*}
\begin{split}
&{_0^CD}_t^{\alpha(t)}y(t)+\sin(t)y^2(t)=g(t),
\quad 0<t\leq 1,\quad 0<\alpha(t)\leq 1,\\
&y(0)=0,
\end{split}
\end{equation*}
where
\[
g(t)=\frac{\Gamma(\frac{9}{2})}{\Gamma\left(\frac{9}{2}
-\alpha(t)\right)}t^{\frac{7}{2}-\alpha(t)}+\sin(t)t^7.
\]

The exact solution of this problem is $y(t)=t^{\frac{7}{2}}$. 
By considering $\alpha(t)=1-0.5\exp(-t)$, we solved this 
problem with different values of $M$. The numerical results are 
displayed in Figure~\ref{fig:1} and Table~\ref{tab:1}. 
In Figure~\ref{fig:1}, the approximate solutions obtained with 
$M=1,~2,~3$, together with the exact solution of this problem, 
are plotted. Furthermore, by considering $M=2,~6,~10$, the absolute errors 
at some selected points are reported in Table~\ref{tab:1}. From these 
results, the convergence of the numerical solutions to the exact 
one can be easily seen. 
\begin{figure}[H]
\includegraphics[scale=1]{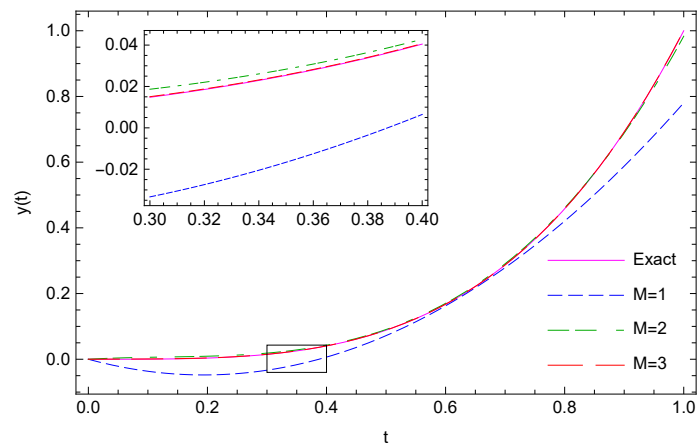}
\caption{(Example~\ref{ex2}) Comparison between the exact solution 
and numerical solutions with $M=1,~2,~3$.}\label{fig:1}
\end{figure}
\vspace{-6pt}
\begin{specialtable}[H]
\setlength{\tabcolsep}{10mm}
\caption{(Example~\ref{ex2}) Absolute errors at some selected 
points with different $M$.}\label{tab:1}
\begin{tabular}{llll}
\toprule
\boldmath{$t$} & \boldmath{$M=2$} & \boldmath{$M=6$} & \boldmath{$M=10$}\\
\midrule
$0.2$ & $5.69\times {10}^{-3}$ & $9.75\times {10}^{-6}$ & $8.06\times {10}^{-7}$ \\
$0.4$ & $2.34\times {10}^{-3}$ & $8.02\times {10}^{-6}$ & $6.34\times {10}^{-7}$ \\
$0.6$ & $2.78\times {10}^{-3}$ & $7.03\times {10}^{-6}$ & $5.53\times {10}^{-7}$ \\
$0.8$ & $2.52\times {10}^{-3}$ & $5.97\times {10}^{-6}$ & $4.59\times {10}^{-7}$ \\
$1.0$ & $1.66\times {10}^{-2}$ & $2.89\times {10}^{-5}$ & $1.95\times {10}^{-6}$ \\ 
\bottomrule
\end{tabular}
\end{specialtable}
\end{Example}


\begin{Example}
\label{ex3}
Consider the following variable-order FDE taken from~\cite{Wu}:
\begin{equation}
\label{6.4}
\begin{split}
&{_0^CD}_t^{\sin(t)}y(t)+y(t)+e^ty(t^5)=g(t),\quad 0<t\leq 1,\\
&y(0)=0,
\end{split}
\end{equation}
where
\begin{equation*}
g(t)=\frac{\Gamma(4)}{\Gamma(4-\sin(t))}t^{3-\sin(t)}
+\frac{\Gamma(3)}{\Gamma(3-\sin(t))}
t^{2-\sin(t)}+e^t(t^{15}+t^{10})+t^3+t^2.
\end{equation*}

The exact solution of this problem is $y(t)=t^3+t^2$. To solve Problem \eqref{6.4}, 
we applied the method with $M=1$ and $M=2$. The numerical solution obtained with $M=1$, 
together with the exact solution are plotted in Figure~\ref{fig:2}. With $M=2$, 
according to the method described in Section~\ref{sec:4}, we set 
$n=\max_{0< t\leq 1}\{\lceil\sin(t)\rceil\}=1$. Therefore, by assuming that
\begin{equation*}
y'(t)=A^TB(t),
\end{equation*}
and using the initial condition, we have:
\begin{equation*}
y(t)=A^TP_t^1B(t),
\quad {_0^CD}_t^{\sin(t)}y(t)=A^TP_t^{1-\sin(t)}B(t),
\quad y(t^5)=A^TP_{t^5}^1B(t^5).
\end{equation*}

By substituting these approximations into \eqref{6.4} 
and using the collocation points $t_0=\frac{1}{4}$, $t_1=\frac{1}{2}$, and $t_2=\frac{3}{4}$, 
we obtain a system of three nonlinear algebraic equations in terms of the elements 
of the vector $A$. By solving the resulting system, one obtains:
\begin{equation*}
a_0=2,\quad a_1=5,
\quad a_2=3.
\end{equation*}

Finally, using these values, we obtain:
\begin{equation*}
y(t)=\left[\begin{array}{ccc}2&5& 3\end{array}\right]\left[
\begin{array}{ccc}
t & 0 & 0 \\
-\frac{t}{4} & \frac{t}{2} & 0 \\
\frac{t}{36} & -\frac{t}{6} & \frac{t}{3}
\end{array}
\right]\left[
\begin{array}{c}
1 \\
t-\frac{1}{2} \\
t^2-t+\frac{1}{6} 
\end{array}
\right]=t^3+t^2,
\end{equation*}
which is the exact solution.
\begin{figure}[H]
\includegraphics[scale=1]{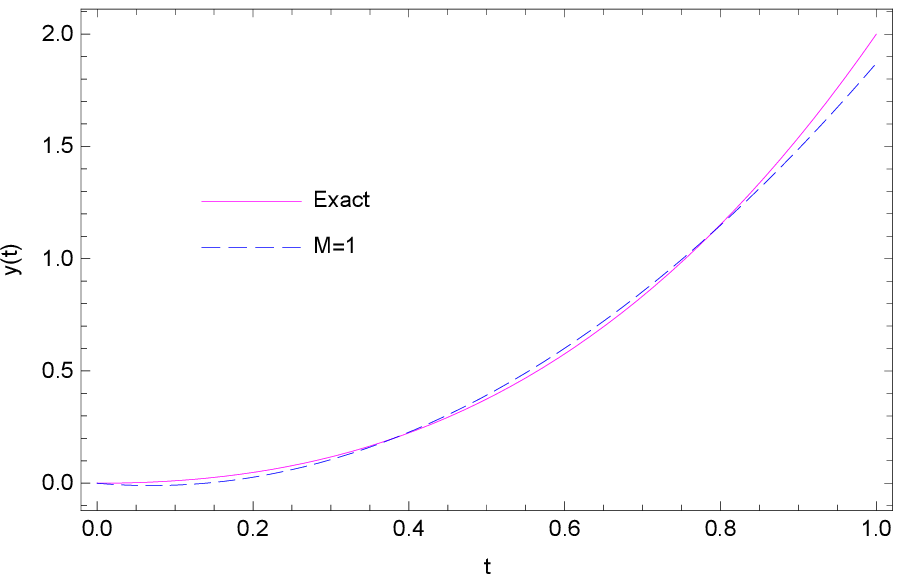}
\caption{(Example~\ref{ex3}) Comparison between the exact solution 
and numerical solution with $M=1$.}\label{fig:2}
\end{figure}
\end{Example}


\begin{Example}
\label{ex4}
Let us now consider the fractional pantograph differential equation: 
\begin{equation}
\label{6.5}
\begin{split}
&{_0^CD}_t^{\alpha(t)}y(t)+y(t)-0.1y(0.2t)=-0.1e^{-0.2t},
\quad 0<t\leq 1,\quad 0<\alpha(t)\leq 1,\\
&y(0)=1.
\end{split}
\end{equation}

The exact solution of this problem, when $\alpha(t)=1$, 
is $y(t)=e^{-t}$. By choosing $M=1$ and $\alpha(t)=1$, we set:
\begin{equation*}
y'(t)=A^TB(t).
\end{equation*}

Then, by considering the initial condition, we have:
\begin{equation*}
y(t)=A^TP_t^1B(t)+1,
\quad y(0.2)=A^TP_{0.2t}^1B(0.2t)+1.
\end{equation*}

By substituting these approximations into \eqref{6.5}, 
and using the collocation points
$$
t_0=\frac{1}{3},\quad t_1=\frac{2}{3},
$$
we obtain
\begin{equation*}
a_0=-0.620328,\quad a_1=0.621053,
\end{equation*}
which gives
\begin{equation*}
y(t)=\left[
\begin{array}{cc}
-0.620328&0.621053
\end{array}\right]\left[
\begin{array}{cc}
t & 0 \\
-\frac{t}{4} & \frac{t}{2}
\end{array}
\right]\left[
\begin{array}{c}
1 \\
t-\frac{1}{2}
\end{array}
\right]+1=0.310526 t^2-0.930854 t+1.
\end{equation*}

This approximate solution and the exact solution to the problem, 
corresponding to $\alpha=1$, are displayed in Figure~\ref{fig:3}. 
By computing the $L^2$-norm of the error for this approximation, we have:
\begin{equation*}
\left\|e^{-t}-(0.310526 t^2-0.930854 t+1)\right\|_2=6.29\times {10}^{-3},
\end{equation*}
which shows that the method gives a high-accuracy approximate solution, 
even with a small number of basis functions. A comparison of the absolute errors, 
obtained by the proposed method with $M=6, 8, 10$ at some selected points, 
with the results proposed in~\cite{Nemati2}, using modified hat functions, 
and those of~\cite{Rahimkhani1}, using Bernoulli wavelets, are reported 
in Table~\ref{tab:2}. From this table, it is seen that our method gave 
more accurate results with a smaller number of basis functions when compared 
to previous methods. Moreover, the approximate solutions obtained with $M=2$ 
and different $\alpha(t)$, along with the exact solution of corresponding 
first-order equation, are given in Figure~\ref{fig:4}. This figure shows 
that the numerical solution is close to the exact solution for 
the case $\alpha(t)=1$ when $\alpha(t)$ is close to one.
\begin{figure}[H]
\includegraphics[scale=1]{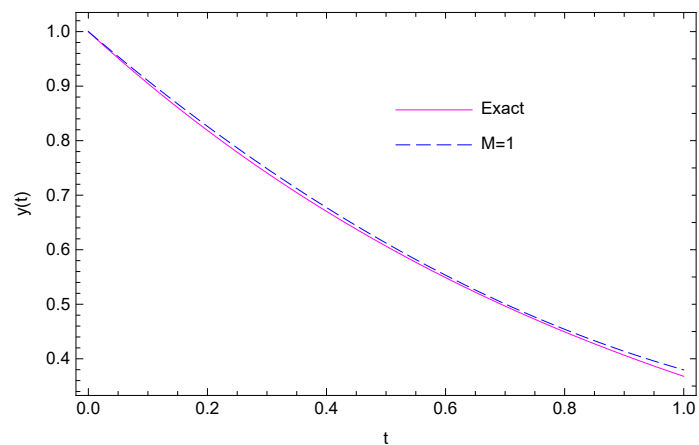}
\caption{(Example~\ref{ex4}) Comparison between exact 
and numerical solutions with $M=1$ and $\alpha(t)=1$.}\label{fig:3}
\end{figure}
\vspace{-6pt}
\begin{specialtable}[H]
\setlength{\tabcolsep}{2.4mm}
\caption{(Example~\ref{ex4}) Comparison of absolute errors, 
at some selected grid points, with $\alpha(t)=1$.}\label{tab:2}
\begin{tabular}{llllllll}
\toprule
&\multicolumn{2}{c}{\textbf{Method of~\cite{Nemati2}}}
&\multicolumn{2}{c}{\textbf{Method of~\cite{Rahimkhani1}}}
&\multicolumn{3}{c}{\textbf{Present Method}}\\
\midrule
\boldmath{$t$} &\multicolumn{2}{l}{\boldmath{$n=64$}} 
& \multicolumn{2}{l}{\boldmath{$k=2,~M=6$}}
&\boldmath{$M=6$} & \boldmath{$M=8$} & \boldmath{$M=10$}\\ \midrule
$2^{-2}$& \multicolumn{2}{l}{$1.18\times {10}^{-9}$} 
&\multicolumn{2}{l}{$1.05 \times {10}^{-8}$} 
&$8.61 \times {10}^{-9}$ & $1.37 \times {10}^{-11}$&$5.56 \times {10}^{-13}$\\
$2^{-3}$& \multicolumn{2}{l}{$5.39\times {10}^{-10}$} 
&\multicolumn{2}{l}{$5.79\times {10}^{-9}$} 
&$1.01 \times {10}^{-8}$ & $1.57 \times {10}^{-11}$&$4.25 \times {10}^{-13}$\\
$2^{-4}$& \multicolumn{2}{l}{$1.17\times {10}^{-9}$} 
&\multicolumn{2}{l}{$2.00\times {10}^{-8}$} &$9.30 \times {10}^{-9}$ 
&$1.59 \times {10}^{-11}$ &$2.42 \times {10}^{-13}$\\
$2^{-5}$& \multicolumn{2}{l}{$5.34\times {10}^{-10}$} 
&\multicolumn{2}{l}{$3.70 \times {10}^{-9}$} & $6.47 
\times {10}^{-9}$ &$1.21 \times {10}^{-11}$ &$1.29 \times {10}^{-13}$\\
$2^{-6}$& \multicolumn{2}{l}{$2.27\times {10}^{-9}$} 
&\multicolumn{2}{l}{$2.03\times {10}^{-8}$} & $3.83 
\times {10}^{-9}$ &$7.58\times {10}^{-12}$ &$6.72 \times {10}^{-14}$\\ \bottomrule
\end{tabular}
\end{specialtable}
\begin{figure}[H]
\includegraphics[scale=1]{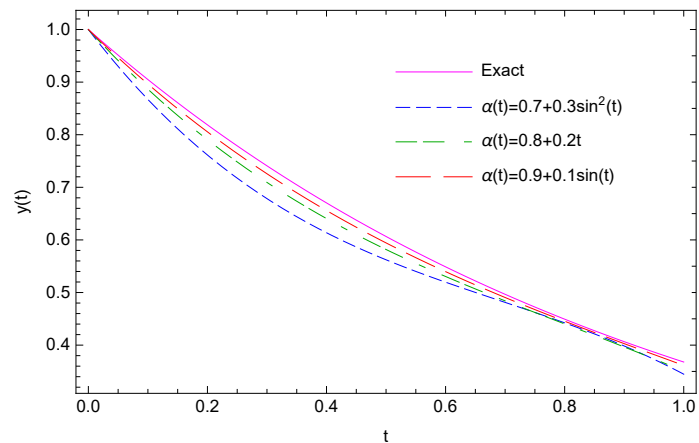}
\caption{(Example~\ref{ex4}) Approximate solutions with $M=2$ 
and different $\alpha(t)$, together with the exact solution 
for $\alpha(t)=1$.}\label{fig:4}
\end{figure}
\end{Example}


\begin{Example}
\label{ex5}
As our last example, consider the following variable-order 
FDE found in~\cite{saber}:
\begin{align*}
{_0^CD}_t^{\alpha(t)}y(t)+3y'(t)-y(t)
=e^t\left(3-\frac{\Gamma(1-\alpha(t),t)}{\Gamma(1-\alpha(t))}\right),
\quad 0<t\leq 1,
\end{align*}
where $\alpha(t)=0.25(1+\cos^2(t))$ and the initial condition is $y(0)=1$. 
The exact solution is $y(t)=e^t$. We solved this problem using the 
proposed method in this paper with different values of $M$. A comparison 
of the absolute errors at some selected points, obtained by our method 
and the one of~\cite{saber}, based on a class of Lagrange polynomials, 
is given in Table~\ref{tab:3}. 
\begin{specialtable}[H]
\setlength{\tabcolsep}{3.4mm}
\caption{(Example~\ref{ex5}) Comparison of the absolute errors 
at some selected grid points.}\label{tab:3}
\begin{tabular}{llllll}
\toprule
&\multicolumn{2}{c}{\textbf{Method of~\cite{saber}}}
&\multicolumn{3}{c}{\textbf{Present Method}}\\
\midrule
\boldmath{$t$} &\boldmath{$M=6$} &\boldmath{$M=10$} &\boldmath{$M=6$} 
& \boldmath{$M=8$} & \boldmath{$M=10$}\\ \midrule
$0.1$& $8.66 \times {10}^{-9}$ & $1.04 \times {10}^{-12}$
&$2.56 \times {10}^{-8}$&$4.12 \times {10}^{-11}$ & $4.40 \times {10}^{-14}$\\
$0.3$& $1.60\times{10}^{-8}$ & $4.57 \times {10}^{-14}$ 
&$2.43 \times {10}^{-8}$ &$3.92 \times {10}^{-11}$& $4.23 \times {10}^{-14}$\\
$0.5$& $2.49\times {10}^{-8}$ &$2.82 \times {10}^{-11}$ 
&$2.44 \times {10}^{-8}$ &$3.93 \times {10}^{-11}$ &$4.24 \times {10}^{-14}$\\
$0.7$& $4.19\times {10}^{-8}$ &$3.12 \times {10}^{-11}$ 
&$2.47 \times {10}^{-8}$ &$3.98 \times {10}^{-11}$&$4.29 \times {10}^{-14}$\\
$0.9$& $5.93\times {10}^{-8}$ &$1.46 \times {10}^{-10}$ &$2.56 \times {10}^{-8}$ 
&$4.14 \times {10}^{-11}$&$4.43\times {10}^{-14}$\\ \bottomrule
\end{tabular}
\end{specialtable}
\end{Example}


\section{Concluding Remarks}
\label{sec:7}

Many researchers have employed fractional differential equations (FDEs) 
in order to model and analyze various scientific phenomena. Typically,
such FDEs do not have known analytical solutions, and approximate 
and numerical approaches have to be applied~\cite{rev01:02}. Here, 
a new numerical method, based on Bernoulli polynomials, was presented 
for solving multiterm variable-order fractional differential equations. 
The operational matrix of variable-order fractional integration for the 
Bernoulli basis functions was introduced, which is a lower triangular 
matrix and helps to reduce the computational effort of the method. 
Our scheme uses this matrix to give some approximations of the unknown 
solution of the problem and its variable-order fractional derivatives 
in terms of the Bernoulli functions. Substituting these approximations 
into the equation and using some collocation points allowed us to reduce 
the problem to a system of nonlinear algebraic equations, 
which greatly simplifies the problem. An error estimate of the method 
was proven, and the applicability of our method was illustrated by solving 
five illustrative examples. The obtained results confirmed the efficiency, 
accuracy, and high performance of our technique, when compared with 
the state-of-the-art numerical schemes available in the literature. 
We emphasize that the accuracy of our method is preserved even when 
the solution of the problem is not infinitely differentiable. This can 
be observed in Example~\ref{ex2}, whose exact solution is $t^{\frac{7}{2}}$.

We used a variable-order definition where the operator 
has no order memory, a so-called type-I operator. As future work, 
it would be interesting to extend the proposed numerical method 
to approximate variable-order derivatives with weak (type-II) 
and strong (type-III) variable-order definitions
\cite{rev03:09,rev03:05}. Other interesting lines of research 
include the stability analysis of the proposed numerical method
and its application to different areas in science and engineering, 
for example in structural mechanics. 


\vspace{6pt}

\authorcontributions{Conceptualization, S.N., P.M.L., and D.F.M.T.; 
methodology, S.N., P.M.L., and D.F.M.T.; software, S.N.; 
validation, S.N., P.M.L., and D.F.M.T.; formal analysis, S.N., P.M.L., and D.F.M.T.; 
investigation, S.N., P.M.L., and D.F.M.T.; 
writing---original draft preparation, S.N., P.M.L., and D.F.M.T.; 
writing---review and editing, S.N., P.M.L., and D.F.M.T.; 
visualization, S.N. All authors have read and agreed 
to the published version of the manuscript.}

\funding{This research was funded by 
Funda\c c\~ao para a Ci\^encia e a Tecnologia 
(FCT, the Portuguese Foundation for Science and Technology)	
through CEMAT, Grant Number UIDB/04621/2020 (P.M.L.),
and CIDMA, Grant Number UIDB/04106/2020 (D.F.M.T.).}

\institutionalreview{Not applicable.} 

\informedconsent{Not applicable.} 

\dataavailability{Not applicable.} 

\acknowledgments{The authors are grateful to the three anonymous reviewers
for several constructive remarks, questions, and suggestions.}

\conflictsofinterest{The authors declare no conflict of interest. 
The funders had no role in the design of the study; in the collection, 
analyses, or interpretation of the data; in the writing of the manuscript; 
nor in the decision to publish the~results.} 


\end{paracol}

\reftitle{References}


\end{document}